\documentclass[a4paper,12pt]{article}

\usepackage[utf8]{inputenc}
\usepackage[T1]{fontenc}
\usepackage[english]{babel}

\usepackage{amssymb,amsthm}
\usepackage{mathtools}

\usepackage{tikz}
\usepackage{float}

\newtheorem{definition}{Definition}[section]

\newtheorem{remark}{Remark}[section]
\newtheorem*{remark*}{Remark}

\newtheorem*{example*}{Example}

\newtheorem{theorem}{Theorem}[section]

\newtheorem{cor}{Corollary}[section]
\newtheorem*{cor*}{Corollary}

\newtheorem{lem}{Lemma}[section]
\newtheorem*{lem*}{Lemma}

\newtheorem*{conj*}{Conjecture}

\newtheorem*{prop*}{Proposition}

\title{The convergence of a gesture recognizer and the shape of a plane gesture}

\author{Lorenzo Luzzi and Paolo Roselli}

\begin{document}
\maketitle

\begin{abstract}
In this work we develop the mathematical framework of !FTL, a new gesture recognition algorithm, published in~\cite{2018-!FTL1}, and we prove its convergence. Such convergence suggests to adopt a notion of shape for smooth gestures as a complex valued function. However, the idea inspiring that notion came to us from Clifford numbers and not from complex numbers. Moreover, the Clifford vector algebra can be used to extend to higher dimensions the notion of ``shape'' of a gesture, while complex numbers are useless to that purpose.
\end{abstract}

\section{Introduction}
A new gesture recognition algorithm, named !FTL by J.L. P\'erez-Medina in~\cite{2018-!FTL1}, has a recognition rate aligned to that of the state-of-the-art \$P recognizer family. Besides, !FTL has proved to be 3 time faster than \$P, thanks to its intrinsic invariance with respect to translation, dilation, and rotation. Indeed, such invariance avoids time consuming  rescaling and normalizing pre-processes. 
In the first part of this article we describe the mathematical framework used to implement !FTL. The notion of shape of a basic gesture is fundamental to !FTL, and recalls that of Lester in~\cite{1996-Lester1}. Then, we will show that !FTL is a discretized version of a limit functional which measures the variation between the  shapes of two plane gestures; thus, extending Lester's notion of shape from triangles to gestures. The proof of convergence will be provided using complex numbers. However,  while preparing this paper, we were faced to some conflicting items:

\begin{itemize}
	\item the ideas inspiring our results come from the geometric  interpretations\footnote{See~\cite{1999-Hestenes} and~\cite{2007-DFM}, for instance.} of the Clifford numbers\footnote{See \cite{1993-Riesz-Lounesto}} , and not from the usual geometry\footnote{See~\cite{1999-Needham}.} of complex numbers;
	\item complex numbers are well known, unlike Clifford numbers;
	\item complex numbers can model plane geometry, while Clifford numbers can model the geometry of a quadratic vector space\footnote{See, for instance,~\cite{1991-Gil-Mur} or~\cite{2013-Meinrenken}.} of any dimension; as a matter of fact, the geometry of Clifford numbers is uniquely determined by the non-degenerate quadratic form defined on the corresponding generating finite-dimensional real vector space; this multidimensional adaptability  allows to extend to higher dimensions the notion of shape of a gesture.
\end{itemize}
As we consider here only gestures in the Euclidean plane, then complex numbers suffice to mimic those   Clifford numbers\footnote{See, for example, the Lounesto's article in~\cite{1993-Riesz-Lounesto}.} in~$\mathcal{C}\ell(2,0)$ (the Clifford vector algebra associated to a two-dimensional Euclidean vector space) needed to state and prove our results. That is why we decided to use complex numbers in this work. Nevertheless,  our results are deeply rooted in the geometric algebra of Clifford numbers. 
For this reason, in the end of the article, we will briefly recall the four dimensional Clifford algebra $\mathcal{C}\ell(2,0)$, and present our results also in that formalism.

\section{Preliminary notions}

A flat surface can be mathematically modeled by the affine\footnote{See, for instance~\cite{2014-Borceux}.} Euclidean plane~$\mathcal{E}$. Thus, the tracing of a single smooth stroke on a flat surface, can be modelled by a function $G:[0,1]\to \mathcal{E}$, such that $G(t)=O+\vec{g}(t)$, where $O\in \mathcal{E}$ is an arbitrary reference point in the plane, $[0,1]=\{t\in\mathbb{R} \ :\ 0\le t\le 1\}$, and 
\[
\vec{g}:  [0,1] \to \mathbb{R}^2 
\]
is a smooth\footnote{That is, twice differentiable on the interval $[0,1]$, with continuous second derivatives.} regular\footnote{That is, whose derivative never vanishes.} vector valued function.
Being the reference point $O$ arbitrary, we can always consider it as the starting point of the gesture, that is $\vec{g}(0)=(0,0)$. In this sense, a gesture is completely determined by the vector-valued function~$\vec{g}$. That is why we give the following definition.

\begin{definition}
A \textbf{plane gesture} is a function $\vec{g}:  [0,1]  \to \mathbb{R}^2 $ which is two times continuously differentiable, and whose derivative $\vec{g}\ '(t)\in \mathbb{R}^2$ is a vector that never vanishes. Briefly $\vec{g}\in C^2\big([0,1];\mathbb{R}^2\big)$, and  $\vec{g}\ '(t)\ne \vec{0}=(0,0)$, for each $t\in[0,1]$.
\end{definition}

When tracing a gesture on a physical device, only a finite number of points are sampled from the input device. A sampled gesture is a finite sequence of  points with timestamps; we describe it mathematically as follows.

\begin{definition}\label{def:n-sample}
A regular \textbf{$n$-sample} of a plane gesture $\vec{g}$ is a sequence of~$n+1$ vectors 
$\vec{g}_0,\dots , \vec{g}_k, \dots ,\vec{g}_n$, where 
\begin{itemize}
	\item $\vec{g}_k= \vec{g}(t_k)$,
	\item $0=t_0< \cdots < t_k< t_{k+1}< \cdots < t_n=1$,
	\item $\Delta\vec{g}_k= \vec{g}(t_k)-\vec{g}(t_{k-1})\ne \vec{0}=(0,0)$, for every $k=1,\dots , n$.
\end{itemize}
\end{definition}

 The following notion of basic gesture is based on the idea that a shape can arise from at least two consecutive sampled points\footnote{Besides the unavoidable  starting point.} of a gesture; that is, two vectors in $\mathbb{R}^2$ (see also Remark~\ref{rem:basic gesture = ordered triangle}).

\begin{definition}
A plane \textbf{basic gesture} is an ordered pair~$(\vec{v}_1,\vec{v}_2)$ of non-zero vectors $\vec{v}_1$, $\vec{v}_2$ in~$\mathbb{R}^2$. 
\end{definition}

\begin{remark}\label{rem:basic gesture = ordered triangle}
A basic gesture can be thought as a particular $2$-sample of a plane gesture  tracing a triangle. More precisely,
we can consider the two vectors of a basic gesture~$(\vec{v}_1,\vec{v}_2)$, as a  pair of consecutive arrows joining three sampled consecutive points 
\begin{center}
$G_0=O$, $G_1=O+\vec{v}_1$, and $G_2=O+\vec{v}_2$,
\end{center}
of a plane gesture $\vec{g}$ tracing the affine plane ordered triangle $G_0G_1G_2$; where $\vec{v}_1=\Delta\vec{g}_1=\vec{g}(t_1)-\vec{g}(t_0)\ne \vec{0}$, $\vec{v}_2=\Delta\vec{g}_2=\vec{g}(t_2)-\vec{g}(t_1)\ne \vec{0}$, and $0=t_0<t_1<t_2=1$. Thus, those three points are the vertexes of a (ordered) triangle, eventually degenerate, whose third oriented side can be traced by the vector $-(\vec{v_1}+\vec{v}_2)$. 

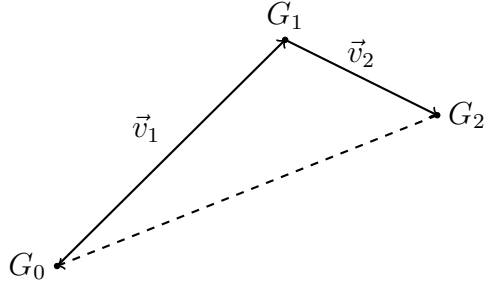
\begin{figure}[H]
\begin{center}
\begin{tikzpicture}[xscale=1,yscale=1]
\draw [thick,->] (0,0)--(3,3);
\draw [thick,->] (3,3)--(5,2);
\draw [thick,dashed,->] (5,2)--(0,0);
\draw [fill] (0,0) circle (1pt);
\draw [fill] (3,3) circle (1pt);
\draw [fill] (5,2) circle (1pt);
\node [above left] at (1.5,1.5) {$\vec{v}_1$};
\node [above] at (4,2.5) {$\vec{v}_2$};
\node [left] at (0,0) {$G_0$};
\node [above] at (3,3) {$G_1$};
\node [right] at (5,2) {$G_2$};
\end{tikzpicture}
\end{center}
\caption{A basic gesture $(\vec{v}_1,\vec{v}_2)$ tracing an ordered affine triangle $G_0G_1G_2$.}\label{fig:ordered triangle}
\end{figure}
One can note that, following the foregoing procedure, a single generic (unordered) triangle can be traced by six basic gestures, possibly different.
\end{remark}

\section{The complex number point of view}

We recall the well known one-to-one correspondence between vectors in $\mathbb{R}^2$ and complex numbers
\begin{equation}
\mathbb{R}^2\ni \vec{v}=(x,y) \longleftrightarrow x+\mathbf{i}y=\mathbf{v}\in\mathbb{C}\ ,
\label{eq:vector<->complex}
\end{equation}
where $x,y\in\mathbb{R}$, and $\mathbf{i}=\sqrt{-1}$, that is $\mathbf{i}^2=-1$. The commutative product between two complex numbers $\mathbf{u}=r+\mathbf{i}s$ and $\mathbf{v}=x+\mathbf{i}y$ is the complex number 
\[
\mathbf{u}\mathbf{v}
=
(rx-sy)+\mathbf{i}(ry+sx)\ .
\]
Thus, if $\mathbf{v}\ne 0$, the quotient between $\mathbf{u}$ and $\mathbf{v}$ is the complex number 
\[
\frac{\mathbf{u}}{\mathbf{v}}=\frac{rx+sy}{x^2+y^2}-\mathbf{i}\frac{ry-sx}{x^2+y^2}\ .
\]

\subsection{The Local Shape Distance}

\begin{definition}
The \textbf{complex shape} of a basic gesture $(\vec{v}_1,\vec{v}_2)$ is the complex number
\[
\displaystyle \frac{\mathbf{v}_1}{\mathbf{v}_2}\ ,
\]
where $\mathbf{v}_1$, $\mathbf{v}_2$ are the two complex numbers corresponding to vectors $\vec{v}_1$, $\vec{v}_2$, respectively, according to the correspondence~(\ref{eq:vector<->complex}).
\end{definition}

\begin{remark*}
The foregoing definition is rooted in Lester's article~\cite{1996-Lester1}, where it is shown that two ordered triangles are similar if and only if the basic gestures generating them (in the sense of Remark~\ref{rem:basic gesture = ordered triangle}) have the same complex shape. 
\end{remark*}

\begin{definition}
The \textbf{Local Shape Distance} between two basic gestures $(\vec{u}_1,\vec{u}_2)$ and $(\vec{v}_1,\vec{v}_2)$ is the non-negative real number
\[
\displaystyle 
LSD\big((\vec{u}_1,\vec{u}_2),(\vec{v}_1,\vec{v}_2)\big)=
\left| \frac{\mathbf{u}_1}{\mathbf{u}_2}- \frac{\mathbf{v}_1}{\mathbf{v}_2}\right|_{\mathbb{C}}\ ,
\]
where $\mathbf{u}_i$, $\mathbf{v}_i$ are the complex numbers corresponding to vectors $\vec{u}_i$, $\vec{v}_i$, respectively, according to the correspondence~(\ref{eq:vector<->complex}). We recall that  
\begin{center}
$
\big|\mathbf{u}-\mathbf{v}\big|_\mathbb{C}
=
\sqrt{(r-x)^2+(s-y)^2}
=
\big|
\vec{u}-\vec{v}
\big|_{\mathbb{R}^2}
$,
\end{center}
where $\mathbf{u}=r+\mathbf{i}s\in \mathbb{C}$, $\mathbf{v}=x+\mathbf{i}y\in\mathbb{C}$, $\vec{u}=(r,s)\in\mathbb{R}^2$, and $\vec{v}=(x,y)\in\mathbb{R}^2$, according to the correspondence~(\ref{eq:vector<->complex}). Thus, the Local Shape Distance is simply the distance between the numbers representing the complex shapes of two basic gestures, according to~\cite{1996-Lester1}. 
\end{definition}

\subsection{The !FTL algorithm}


\begin{definition}\label{def:!FTL}
Given the $n$-samples of two plane gestures $\vec{f}$ and $\vec{g}$, the following \textbf{!FTL algorithm} gives a measure of their dissimilarity based on the Local Shape Distance between the complex shapes of basic gestures produced by consecutive pairs of vectors taken from the samples. More precisely, given two isochronous\footnote{That is, they are sampled at the same timestamps 	$0=t_0< \cdots < t_n=1$.} $n$-samples\footnote{See Definition~\ref{def:n-sample}.} 
\begin{center}
$
\vec{f}_0, \ \dots\ , \vec{f}_n \ \ ,\ \
\vec{g}_0, \ \dots\ , \vec{g}_n
$
\end{center}
of the plane gestures $\vec{f}(t)=\big(r(t),s(t)\big)$, $\vec{g}(t)=\big(x(t),y(t)\big)\in\mathbb{R}^2$, respectively, where $r,s,x,y \in C^2\big([0,1];\mathbb{R}\big)$, then 
\begin{align*}
!FTL(\vec{f}_0, \ \dots\ , \vec{f}_n\ ,\ 
\vec{g}_0, \ \dots\ , \vec{g}_n)
& =
\sum_{k=1}^{n-1}
LSD\big((\Delta\vec{f}_k,\Delta\vec{f}_{k+1}),(\Delta\vec{g}_k,\Delta\vec{g}_{k+1})\big)\\
& =
\sum_{k=1}^{n-1}
\left|\frac{\Delta\mathbf{f}_k}{\Delta\mathbf{f}_{k+1}}-\frac{\Delta\mathbf{g}_k}{\Delta\mathbf{g}_{k+1}}\right|_\mathbb{C}\ .
\end{align*}
where $\Delta\mathbf{f}_k=\mathbf{f}(t_k)-\mathbf{f}(t_{k-1})$, $\Delta\mathbf{g}_k=\mathbf{g}(t_k)-\mathbf{g}(t_{k-1})$, $\mathbf{f}(t)=r(t)+\mathbf{i}s(t)$, and $\mathbf{g}(t)=x(t)+\mathbf{i}y(t)$ are complex numbers.
\end{definition}

Some natural questions about !FTL arise,
\begin{itemize}
	\item if gesture $\vec{g}$ is translated, does the value of !FTL change? 
	\item if gesture $\vec{g}$ is uniformly scaled, does the value of !FTL change? 
	\item if gesture $\vec{g}$ is rotated, does the value of !FTL change? 
	\item if one increases the number of the sampled points of the two gestures, does !FTL has a limit value?
	\item if such a limit value exists, can we find a closed formula to express it?
\end{itemize}

In what follows, we will see that !FTL is invariant with respect to translations, scaling and rotations; moreover, under certain hypothesis, to increase the number of sampled points improves the measure of dissimilarity, which corresponds to a well definite number explicitly expressed as a Riemann integral. 

\subsubsection{Invariance properties of LSD}

In this section, we will show the invariance properties of the complex shape $\displaystyle \frac{\Delta\mathbf{g}_k}{\Delta\mathbf{g}_{k+1}}$ of a basic gesture $(\Delta\vec{g}_k,\Delta\vec{g}_{k+1})$ coming from the sample of a plane gesture $\vec{g}$.
As a matter of fact, 
let $\vec{g}_0, \ \dots\ , \vec{g}_n$ be an $n$-sample of a plane gesture $\vec{g}(t)=\big(x(t),y(t)\big)$; then
\begin{enumerate}
	\item for each vector $\vec{v}\in\mathbb{R}^2$,
		\begin{itemize}
			\item $\vec{p}(t)=\vec{g}(t)+\vec{v}$ is a plane gesture, 
			\item $\Delta \vec{p}_k=\Delta \vec{g}_k$, for all $k=1,\dots , n$,
		\end{itemize}
	and then $\displaystyle \frac{\Delta \mathbf{p}_k}{\Delta \mathbf{p}_{k+1}}=\frac{\Delta \mathbf{g}_k}{\Delta \mathbf{g}_{k+1}}$.
	\item for each number $\lambda\in\mathbb{R}$ $(\lambda \ne 0)$, 
		\begin{itemize}
			\item $\vec{l}(t)=\lambda \vec{g}(t)$ is a plane gesture, 
			\item $\Delta \vec{l}_k=\lambda \Delta \vec{g}_k$, for all $k=1,\dots , n$,
		\end{itemize}
	and then 
	$\displaystyle 
	\frac{\Delta \mathbf{l}_k}{\Delta \mathbf{l}_{k+1}}
	=
	\frac{\lambda \Delta \mathbf{g}_k}{\lambda\Delta \mathbf{g}_{k+1}}
	=
	\frac{\Delta \mathbf{g}_k}{\Delta \mathbf{g}_{k+1}}$.
	\item for each number $\theta\in\mathbb{R}$, 
		\begin{itemize}
			\item $\vec{q}(t)=\big( x(t)\cos\theta - y(t)\sin\theta \ ,\ y(t)\cos\theta + x(t)\sin\theta \big)$ is a plane gesture, 
			\item $\mathbf{q}(t)= (\cos\theta +\mathbf{i}\sin\theta)(x(t)+\mathbf{i}y(t))=e^{\mathbf{i}\theta}\mathbf{g}(t)$
			\item $\Delta \mathbf{q}_k= e^{\mathbf{i}\theta}\Delta\mathbf{g}_k$, for all $k=1,\dots , n$,
		\end{itemize}
	and then 
	$
	\displaystyle 
	\frac{\Delta \mathbf{q}_k}{\Delta \mathbf{q}_{k+1}}
	=
	\frac{e^{\mathbf{i}\theta}\Delta \mathbf{g}_k}{e^{\mathbf{i}\theta}\Delta \mathbf{g}_{k+1}}
	=
	\frac{\Delta \mathbf{g}_k}{\Delta \mathbf{g}_{k+1}}
	$.
\end{enumerate}
The invariance properties of !FTL follows  then from those of the complex shape, thanks to Definition~\ref{def:!FTL}.

\subsubsection{Convergence of !FTL through complex numbers}


\begin{theorem}\label{thm:convergence uniform !FTL}
Given two plane gestures $\vec{f}(t)=\big(r(t),s(t)\big)\in\mathbb{R}^2$, and $\vec{g}(t)=\big(x(t),y(t)\big)\in\mathbb{R}^2$, then
\[
\lim_{n\to \infty}
!FTL(\vec{f}_0, \ \dots\ , \vec{f}_n\ ,\ 
\vec{g}_0, \ \dots\ , \vec{g}_n)
=
\int_0^1
\left|
\frac{\mathbf{f}''(t)}{\mathbf{f}'(t)}-\frac{\mathbf{g}''(t)}{\mathbf{g}'(t)}
\right|_{\mathbb{C}}\ dt
\]
where 
\begin{itemize}
	\item $\vec{f}_0, \ \dots\ , \vec{f}_n$, $\vec{g}_0, \ \dots\ , \vec{g}_n$ are isochronous $n$-samples of $\vec{f}$ and $\vec{g}$, respectively, such that $ t_k = \frac{k}{n}$, for all $k=0,1,\dots n$,
	\item $\mathbf{f}(t)=r(t)+\mathbf{i}s(t)$, and $\mathbf{g}(t)=x(t)+\mathbf{i}y(t)$.
\end{itemize}
\end{theorem}

{\it A proof.
By hypothesis, the Riemann integral $\displaystyle \int_0^1
\left|
\frac{\mathbf{f}''(t)}{\mathbf{f}'(t)}-\frac{\mathbf{g}''(t)}{\mathbf{g}'(t)}
\right|_{\mathbb{C}}\ dt$ exists; this implies that for every $\epsilon>0$ there exists $N_\epsilon \in\mathbb{N}$ such that 
\[
\left|
\sum_{k=1}^n
\left|\frac{\mathbf{f}''(\xi_k)}{\mathbf{f}'(\xi_k)}-\frac{\mathbf{g}''(\xi_k)}{\mathbf{g}'(\xi_k)}\right|_{\mathbb{C}}\frac{1}{n} 
-
\int_0^1
\left|
\frac{\mathbf{f}''(t)}{\mathbf{f}'(t)}-\frac{\mathbf{g}''(t)}{\mathbf{g}'(t)}
\right|_{\mathbb{C}}\ dt
\right|
< \epsilon\ ,
\] 
provided $n > N_\epsilon$, and $\xi_k\in \left[\frac{k-1}{n},\frac{k}{n}\right]$, with $k=1,\dots , n$.

Notice that, to evaluate each shape $\displaystyle \frac{\Delta\mathbf{g}_k}{\Delta\mathbf{g}_{k+1}}$, the extremities of two adjacent intervals are needed. In particular, we can write
\begin{equation}
\label{eq:double sum 1}
\sum_{k=1}^{2m-1}
\frac{\Delta\mathbf{g}_k}{\Delta\mathbf{g}_{k+1}}
=
\sum_{h=1}^{m}
\frac{\Delta\mathbf{g}_{2h-1}}{\Delta\mathbf{g}_{2h}}
+
\sum_{h=1}^{m-1}
\frac{\Delta\mathbf{g}_{2h}}{\Delta\mathbf{g}_{2h+1}}\ ,
\end{equation}
when $n=2m$ is even\footnote{A similar expression holds when $n$ is odd.}.
Thus, to estimate the difference between local shape distances and terms of a Riemann sum, we have to consider  the latter on couples of adjacent intervals. In order to simplify notations, we will consider in the following only the case $n=2m$ ($n$ even). However, our arguments can be applied similarly to the case: $n$ odd.
If $n> 2N_\epsilon$, then the integral can be estimated both by
\begin{itemize}
	\item \[
\left|
\sum_{h=1}^m
\left|\frac{\mathbf{f}''(\xi^e_h)}{\mathbf{f}'(\xi^e_h)}-\frac{\mathbf{g}''(\xi^e_h)}{\mathbf{g}'(\xi^e_h)}\right|_{\mathbb{C}}\frac{1}{n} 
-
\frac{1}{2}
\int_0^1
\left|
\frac{\mathbf{f}''(t)}{\mathbf{f}'(t)}-\frac{\mathbf{g}''(t)}{\mathbf{g}'(t)}
\right|_{\mathbb{C}}\ dt
\right|
< \frac{\epsilon}{2}\ ,
\] 
where $\xi^e_h\in \left[\frac{2(h-1)}{n},\frac{2h}{n}\right]$, with $h=1,\dots , m$, and 
	\item \[
\left|
\sum_{h=1}^{m-1}
\left|\frac{\mathbf{f}''(\xi^o_h)}{\mathbf{f}'(\xi^o_h)}-\frac{\mathbf{g}''(\xi^o_h)}{\mathbf{g}'(\xi^o_h)}\right|_{\mathbb{C}}\frac{1}{n} 
-
\frac{1}{2}
\int_0^1
\left|
\frac{\mathbf{f}''(t)}{\mathbf{f}'(t)}-\frac{\mathbf{g}''(t)}{\mathbf{g}'(t)}
\right|_{\mathbb{C}}\ dt
\right|
< \frac{\epsilon}{2}\ ,
\] 
where $\xi^o_h\in \left[\frac{2h-1}{n},\frac{2h+1}{n}\right]$, with $h=1,\dots , m$.
\end{itemize}

Then, to obtain the thesis, it suffices to see how to estimate the following quantity,
\begin{align*}
&
\left|
\frac{\Delta\mathbf{f}_{2h}}{\Delta\mathbf{f}_{2h+1}}
-
\frac{\Delta\mathbf{g}_{2h}}{\Delta\mathbf{g}_{2h+1}}
-
\left(\frac{\mathbf{g}''(\xi^e_h)}{\mathbf{g}'(\xi^e_h)}-\frac{\mathbf{f}''(\xi^e_h)}{\mathbf{f}'(\xi^e_h)}\right)\frac{1}{n} \right|_{\mathbb{C}}=\\
=& 
\left|
\left(
\frac{\Delta\mathbf{f}_{2h}}{\Delta\mathbf{f}_{2h+1}}-1+\frac{\mathbf{f}''(\xi^e_h)}{\mathbf{f}'(\xi^e_h)}\frac{1}{n}\right)
+
\left( 1 -
\frac{\Delta\mathbf{g}_{2h}}{\Delta\mathbf{g}_{2h+1}}- \frac{\mathbf{g}''(\xi^e_h)}{\mathbf{g}'(\xi^e_h)}\frac{1}{n}\right)
\right|_{\mathbb{C}}
\ ,
\end{align*}
for each $h=1,\dots , m$. In particular\footnote{A similar argument can be applied for the function $\mathbf{g}$.}, we can observe that, assuming $\delta=\frac{1}{n}$, then 
\begin{align*}
1 -
\frac{\Delta\mathbf{g}_{2h}}{\Delta\mathbf{g}_{2h+1}}
& =
1-
\frac{\mathbf{g}(t_{2h})-\mathbf{g}(t_{2h}-\delta)}{\mathbf{g}(t_{2h}+\delta)-\mathbf{g}(t_{2h})}
=
\frac{\mathbf{g}(t_{2h}+\delta)-2\mathbf{g}(t_{2h})+\mathbf{g}(t_{2h}-\delta)}{\mathbf{g}(t_{2h}+\delta)-\mathbf{g}(t_{2h})}=\\
& =
\frac{
\frac{\mathbf{g}(t_{2h}+\delta)-2\mathbf{g}(t_{2h})+\mathbf{g}(t_{2h}-\delta)}{\delta^2}}{\frac{\mathbf{g}(t_{2h}+\delta)-\mathbf{g}(t_{2h})}{\delta}}\delta\ .
\end{align*}
By hypothesis, the function $\mathbf{g}$ is twice differentiable and $\mathbf{g}'\ne 0$, thus we have that, for every $t\in [0,1]$
\[
\lim_{\delta \to 0}
\frac{
\frac{\mathbf{g}(t+\delta)-2\mathbf{g}(t)+\mathbf{g}(t-\delta)}{\delta^2}}{\frac{\mathbf{g}(t+\delta)-\mathbf{g}(t)}{\delta}}
=
\frac{\mathbf{g}''(t)}{\mathbf{g}'(t)}\ ,
\]
as the limit of a quotient is the quotient of the limits, provided the limit of the denominator is not zero. So, we have that, for every $\epsilon >0$ there exists $\delta_\epsilon$, such that if $\delta<\delta_\epsilon$, then  
\begin{align*}
\left|
1
-
\frac{\Delta\mathbf{g}_{2h}}{\Delta\mathbf{g}_{2h+1}}
-
\frac{\mathbf{g}''(\xi_k)}{\mathbf{g}'(\xi_k)}\delta
\right|_{\mathbb{C}}
=
\left|
\frac{
\frac{\mathbf{g}(t_{2h}+\delta)-2\mathbf{g}(t_{2h})+\mathbf{g}(t_{2h}-\delta)}{\delta^2}}{\frac{\mathbf{g}(t_{2h}+\delta)-\mathbf{g}(t_{2h})}{\delta}}
-
\frac{\mathbf{g}''(\xi^e_{h})}{\mathbf{g}'(\xi^e_h)}
\right|_{\mathbb{C}}\delta 
< \epsilon \delta\ ,
\end{align*}
and this prove the thesis, provided $\delta <\min\{\delta_\epsilon, \frac{1}{2 N_\epsilon}\}$.\ $\square$
}

The foregoing proof can also be used to prove other results, such as  the following one.
\begin{cor}\label{cor:complex shape}
Given a plane gesture $\vec{g}$, then\footnote{We adopt the same notations of  Theorem~\ref{thm:convergence uniform !FTL}.}
\begin{equation}
\label{eq:twice complex shape}
\lim_{n \to \infty}
\sum_{k=1}^{n-1}
\frac{\Delta\mathbf{g}_k}{\Delta\mathbf{g}_{k+1}}
=
2
-
\int_0^1
\frac{\mathbf{g}''(t)}{\mathbf{g}'(t)}
\ dt\ \in\mathbb{C}.
\end{equation}
\end{cor}
A more general proof of the foregoing result will be given with Theorem~\ref{thm:convergence to shape}.
Corollary~\ref{cor:complex shape} makes then reasonable to give the following definition.

\begin{definition}
The \textbf{Complex Shape} of a plane gesture $\vec{g}(t)=\big(x(t),y(t)\big)\in\mathbb{R}^2$, $x,y\in C^2\big([0,1];\mathbb{R}\big)$, is the following complex valued function
\[
1
-
\frac{\mathbf{g}''(t)}{2\mathbf{g}'(t)}
\ ,
\]
where $\mathbf{g}(t)=x(t)+\mathbf{i}y(t)$ and $t\in[0,1]$.
\end{definition}

\begin{remark*}
We decided to scale $(\ref{eq:twice complex shape})$ in half so that the complex shape of a rectilinear gesture would be $1$, regardless of whether it is considered as ``basic'' or not. Indeed, $(\ref{eq:twice complex shape})$ is the double of the complex shape simply because of a kind of double counting of intervals in relation~$(\ref{eq:double sum 1})$. 
\end{remark*}

\begin{example*}
The complex shape of the circled plane gesture
\[
\vec{g}(t)=\Big(x_0 + r \cos\big(2\pi (t-\phi)\big)\ ,\  y_0 + r \sin\big(2\pi (t-\phi)\big)\Big)\ ,
\]
is the constant value $1-\pi\mathbf{i}$. Notice that it is independent from the radius $r\in\mathbb{R}^+$, the center $(x_0,y_0)\in\mathbb{R}^2$, and the phase $\phi\in\mathbb{R}$, thanks to the invariant properties of the complex shape of a basic gesture.
\end{example*}

\subsubsection{The case of non-uniformly spaced timestamps}

The !FTL algorithm is suited for uniform $n$-samplings, that is, when~$t_k-t_{k-1}$ is independent of index $k$. However, most of sampling devices are multitasking; this implies that the Central Processing Unit is not always sampling points; so, $t_k-t_{k-1}$ may depend on $k$. In this situation, it is reasonable to explore the {\bf weighted complex shape}
\[ 
\frac{t_{k+1}-t_k}{t_k-t_{k-1}}\frac{\Delta \mathbf{g}_k}{\Delta\mathbf{g}_{k+1}} \in \mathbb{C}\ ,
\]
 of a basic gesture 
\[
(\Delta\vec{g}_k,\Delta\vec{g}_{k+1})
=
\big(\vec{g}(t_k)-\vec{g}(t_{k-1}) , \vec{g}(t_{k+1})-\vec{g}(t_{k})\big)\ ,
\]
taken from the $n$-sample of a plane gesture $\vec{g}$ (without assuming that the timestamps are uniformly spaced). Of course, such weighted complex shape coincides with the complex shape when $t_{k+1}-t_k=t_{k}-t_{k-1}$. Moreover, we will show that, such weighted complex shape is still convergent to the same value of~$(\ref{eq:twice complex shape})$.

\begin{lem}\label{lem:second divided difference}
Given a plane gesture $\vec{g}(t)$ then, for each $t\in (0,1)$, we have that
\[
\lim_{
\begin{array}{c}
	\scriptstyle \tau_0 \to t \ , \ \tau_1 \to t \ , \ \tau_2 \to t \\
	\scriptstyle \tau_0 \ne \tau_1 \ , \ \tau_1 \ne \tau_2\ , \ \tau_2 \ne \tau_0
\end{array}
}
\left(
1-
\frac{\tau_2-\tau_1}{\tau_1 - \tau_0}\ 
\frac{\mathbf{g}(\tau_1)-\mathbf{g}(\tau_0)}{\mathbf{g}(\tau_2)-\mathbf{g}(\tau_1)}\right)
\frac{1}{\tau_2-\tau_0}
=
\frac{1}{2}\frac{\mathbf{g}''(t)}{\mathbf{g}'(t)}
\]
\end{lem}
{\it A proof. After rewriting 
\begin{equation}
\label{eq:double ratio}
\left(\kern-2pt
1-
\frac{\tau_2-\tau_1}{\tau_1 - \tau_0}\ 
\frac{\mathbf{g}(\tau_1)-\mathbf{g}(\tau_0)}{\mathbf{g}(\tau_2)-\mathbf{g}(\tau_1)}\right)
\kern-3pt \frac{1}{\tau_2-\tau_0}
=
\frac{
\frac{\mathbf{g}(\tau_2)-\mathbf{g}(\tau_1)}{\tau_2-\tau_1}-\frac{\mathbf{g}(\tau_1)-\mathbf{g}(\tau_0)}{\tau_1-\tau_0}
}
{\tau_2-\tau_0}
\frac{1}{\frac{\mathbf{g}(\tau_2)-\mathbf{g}(\tau_1)}{\tau_2 - \tau_1}}
\end{equation}
we notice that
\begin{equation}
\label{eq:second div diff}
\frac{
\frac{\mathbf{g}(\tau_2)-\mathbf{g}(\tau_1)}{\tau_2-\tau_1}-\frac{\mathbf{g}(\tau_1)-\mathbf{g}(\tau_0)}{\tau_1-\tau_0}
}
{\tau_2-\tau_0}
\end{equation}
is the second divided difference\footnote{See \cite{2016-Burden-Faires} at page 123.} of the complex valued function $\mathbf{g}$ at points $\tau_0,\tau_1$ and $\tau_2$.
Being the function twice continously differentiable, it suffices to apply the Mean Value Theorem for divided differences\footnote{See Theorem 2.10 in~\cite{1998-Sahoo-Riedel}, at page 60.} to real and imaginary parts of~$\mathbf{g}$, to obtain that
\[
\lim_{
\begin{array}{c}
	\scriptstyle \tau_0 \to t \ , \ \tau_1 \to t \ , \ \tau_2 \to t \\
	\scriptstyle \tau_0 < \tau_1 < \tau_2\ 
\end{array}
}
\frac{
\frac{\mathbf{g}(\tau_2)-\mathbf{g}(\tau_1)}{\tau_2-\tau_1}-\frac{\mathbf{g}(\tau_1)-\mathbf{g}(\tau_0)}{\tau_1-\tau_0}
}
{\tau_2-\tau_0}
=
\frac{\mathbf{g}''(t)}{2}\ .
\]
Notice that we can always assume condition $\tau_0 < \tau_1 < \tau_2$; as a matter of fact, the second divided difference~$(\ref{eq:second div diff})$ is symmetric with respect points $\tau_0, \tau_1$ and $\tau_2$.
As the limit of quotient~$(\ref{eq:double ratio})$ is the quotient of the limits, provided the limit of the denominator is not zero, one obtains the thesis. \ $\square$}

\begin{theorem}\label{thm:convergence to shape}
Given a plane gesture $\vec{g}(t)=\big(x(t),y(t)\big)\in\mathbb{R}^2$, then
\[
\lim_{\delta\to 0^+}
\sum_{k=1}^{n-1}
\frac{t_{k+1}-t_k}{t_k-t_{k-1}}
\frac{\Delta\mathbf{g}_k}{\Delta\mathbf{g}_{k+1}}
=
2
-
\int_0^1
\frac{\mathbf{g}''(t)}{\mathbf{g}'(t)}
\ dt
\ \in\mathbb{C}\ ,
\]
where $0=t_0<\cdots <t_{k-1}< t_k<\cdots < t_n=1$, and  $\displaystyle \delta = \max_{1\le k\le n}\{t_k-t_{k-1}\}$.
\end{theorem}

{\it A proof.
 By hypothesis, the complex valued Riemann integral $\displaystyle \int_0^1 \frac{\mathbf{g}''(t)}{\mathbf{g}'(t)}\ dt$ exists; this implies that for every $\epsilon>0$ there exists $\delta_\epsilon>0$ such that 
\[
\left|
\int_0^1
\frac{\mathbf{g}''(t)}{\mathbf{g}'(t)}\ dt
-
\sum_{k=1}^n
\frac{\mathbf{g}''(\xi_k)}{\mathbf{g}'(\xi_k)} (t_k-t_{k-1}) 
\right|_{\mathbb{C}}
< \epsilon\ ,
\] 
provided the partition 
\begin{center}
$0=t_0<\cdots <t_{k-1}< t_k<\cdots < t_n=1$
\end{center}
is such that $t_{k}-t_{k-1}< \delta_\epsilon$, and $\xi_k\in[t_{k-1},t_k]$ for each $k=1\dots , n$.\\
Notice that, to evaluate each shape $\displaystyle \frac{\Delta\mathbf{g}_k}{\Delta\mathbf{g}_{k+1}}$, the extremities of two adjacent intervals are needed. This implies that
\begin{equation}
\label{eq:double sum 2}
\sum_{k=1}^{2m-1}
{\scriptstyle \frac{t_{k+1}-t_k}{t_k-t_{k-1}}}
\frac{\Delta\mathbf{g}_k}{\Delta\mathbf{g}_{k+1}}
=
\sum_{h=1}^{m}
{\scriptstyle \frac{t_{2h}-t_{2h-1}}{t_{2h-1}-t_{2(h-1)}}}
\frac{\Delta\mathbf{g}_{2h-1}}{\Delta\mathbf{g}_{2h}}
+
\sum_{h=1}^{m-1}
{\scriptstyle \frac{t_{2h+1}-t_{2h}}{t_{2h}-t_{2h-1}}}
\frac{\Delta\mathbf{g}_{2h}}{\Delta\mathbf{g}_{2h+1}}\ ,
\end{equation}
when $n$ is even\footnote{A similar expression old when $n$ is odd.}. 
Thus, to estimate the difference between complex shapes and Riemann sums, we need to consider  the latter on couples of adjacent intervals; one with even-indexed extremities, the other with odd-indexed extremities. 
In order to simplify notations, we will consider in the following only partitions of $[0,1]$ having an even number of points ($n=2m$), that is    
\begin{equation}
\label{eq:even partition}
0=t_0 < \cdots <t_{k-1}< t_k <\cdots < t_{2m}=1\ .
\end{equation}
However, our arguments can be applied similarly to partitions of $[0,1]$ having an odd number of points.
If partition $(\ref{eq:even partition})$ is such that 
\[
\displaystyle 
\max\left\{\max_{1\le h\le m}(t_{2h}-t_{2(h-1)})\ , \  \max_{1\le h\le m}(t_{2h+1}-t_{2h-1})\right\}< \delta_\epsilon\ ,
\] 
then we can estimate the Riemann sum both 
\begin{itemize}
	\item on ``even indexed'' intervals
\[
\left|
\int_0^1
\frac{\mathbf{g}''(t)}{\mathbf{g}'(t)}\ dt
-
\sum_{h=1}^m
\frac{\mathbf{g}''(\xi^e_h)}{\mathbf{g}'(\xi^e_h)} (t_{2h}-t_{2(h-1)}) 
\right|
< \epsilon\ ,
\] 
whatever are $\xi^e_h\in[t_{2(h-1)},t_{2h}]$ when $h=1,\dots, m$, and 
\item on ``odd indexed'' intervals, where a similar estimate is possible
\[
\left|
\int_0^1
\frac{\mathbf{g}''(t)}{\mathbf{g}'(t)}\ dt
-
\frac{\mathbf{g}''(\xi^o_0)}{\mathbf{g}'(\xi^o_0)} {\scriptstyle (t_{1}-t_{0})}
-
\frac{\mathbf{g}''(\xi^o_m)}{\mathbf{g}'(\xi^o_m)} {\scriptstyle (t_{2m}-t_{2m-1})}
-
\sum_{h=1}^{m-1}
\frac{\mathbf{g}''(\xi^o_h)}{\mathbf{g}'(\xi^o_h)} {\scriptstyle (t_{2h+1}-t_{2h-1}) }
\right|
< \epsilon\ ,
\] 
whatever are $\xi^o_h\in[t_{2h-1},t_{2h+1}]$, with $h=1,\dots, m-1$, $\xi^o_0\in [t_0,t_1]$, and $\xi^o_m\in [t_{2m-1},t_{2m}]$. 
\end{itemize}
Now, let us focus on the first term of the right expression in~$(\ref{eq:double sum 2})$. In order to get the thesis, we need to estimate  each term
\begin{align*}
& &
\frac{t_{2h}-t_{2h-1}}{t_{2h-1}-t_{2(h-1)}}
\frac{\Delta\mathbf{g}_{2h-1}}{\Delta\mathbf{g}_{2h}}
-
(t_{2h}-t_{2(h-1)})
+
\frac{1}{2}
\frac{\mathbf{g}''(\xi^e_h)}{\mathbf{g}'(\xi^e_h)} (t_{2h}-t_{2(h-1)}) 
=\\
=
& &
\left(
\frac{t_{2h}-t_{2h-1}}{t_{2h-1}-t_{2(h-1)}}
\frac{\Delta\mathbf{g}_{2h-1}}{\Delta\mathbf{g}_{2h}}
\frac{1}{t_{2h}-t_{2(h-1)}}
-
1
+
\frac{1}{2}
\frac{\mathbf{g}''(\xi^e_h)}{\mathbf{g}'(\xi^e_h)}
\right) (t_{2h}-t_{2(h-1)}) \ .
\end{align*}
If one considers {\rm Lemma~\ref{lem:second divided difference}} with $\tau_0=t_{2(h-1)}$, $\tau_1=t_{2h-1}$, and $\tau_2=t_{2h}$, we have the estimate
\[
\left|
\frac{t_{2h}-t_{2h-1}}{t_{2h-1}-t_{2(h-1)}}
\frac{\Delta\mathbf{g}_{2h-1}}{\Delta\mathbf{g}_{2h}}
\frac{1}{t_{2h}-t_{2(h-1)}}
-
1
+
\frac{1}{2}
\frac{\mathbf{g}''(\xi^e_h)}{\mathbf{g}'(\xi^e_h)}
\right|
<
\frac{\epsilon}{2}\ ,
\]
which is independent from index~$h$, thanks to the uniform continuity of $\displaystyle \frac{\mathbf{g}''}{\mathbf{g}'}$. By applying the same lemma for the odd terms involving $\xi^o_h$, the thesis follows.~$\square$
}

\begin{remark*}
The foregoing theorem provide a new algorithm (let us call\footnote{The letter ``W'' stands for ``weighted''.} it~!WFTL), to measure the dissimilarity between the not-necessarily uniformly sampled isochronous $n$-samples, of two plane gestures~$\vec{f}$ and~$\vec{g}$ 
 \begin{align*}
!WFTL(\vec{f}_0, \ \dots\ , \vec{f}_n\ ,\ 
\vec{g}_0, \ \dots\ , \vec{g}_n)
& =
\sum_{k=1}^{n-1}
{\scriptstyle \frac{t_{k+1}-t_k}{t_k-t_{k-1}}}
LSD\big((\Delta\vec{f}_k,\Delta\vec{f}_{k+1}),(\Delta\vec{g}_k,\Delta\vec{g}_{k+1})\big)\\
& =
\sum_{k=1}^{n-1}
{\scriptstyle \frac{t_{k+1}-t_k}{t_k-t_{k-1}}}
\left|\frac{\Delta\mathbf{f}_k}{\Delta\mathbf{f}_{k+1}}-\frac{\Delta\mathbf{g}_k}{\Delta\mathbf{g}_{k+1}}\right|_\mathbb{C}\ .
\end{align*}
\end{remark*}

We can then claim the following result.

\begin{cor*}
Given two plane gestures $\vec{f}(t)=\big(r(t),s(t)\big)\in\mathbb{R}^2$, and $\vec{g}(t)=\big(x(t),y(t)\big)\in\mathbb{R}^2$, then
\[
\lim_{\delta \to 0^+}
!WFTL(\vec{f}_0, \ \dots\ , \vec{f}_n\ ,\ 
\vec{g}_0, \ \dots\ , \vec{g}_n)
=
\int_0^1
\left|
\frac{\mathbf{f}''(t)}{\mathbf{f}'(t)}-\frac{\mathbf{g}''(t)}{\mathbf{g}'(t)}
\right|_{\mathbb{C}}\ dt
\]
where 
\begin{itemize}
	\item $\vec{f}_0, \ \dots\ , \vec{f}_n$, $\vec{g}_0, \ \dots\ , \vec{g}_n$ are isochronous~$n$-samples of $\vec{f}$ and $\vec{g}$, respectively,
	\item $0=t_0<\cdots <t_{k-1}< t_k<\cdots < t_n=1$, $\displaystyle \delta = \max_{1\le k\le n}\{t_k-t_{k-1}\}$,
	\item $\mathbf{f}(t)=r(t)+\mathbf{i}s(t)$, and $\mathbf{g}(t)=x(t)+\mathbf{i}y(t)$.
\end{itemize}
\end{cor*}

\section{The Clifford number point of view}
As we have seen, in order to define the complex shape of a basic gesture, we had to consider the components of a vector $\vec{v}=(x,y)\in\mathbb{R}^2$ as the real and imaginary parts of the complex number $\mathbf{v}=x+\mathbf{i}y$. This twisted construction allows us to use the quotient between complex numbers to encode the concept of ``shape'' as a complex number. However, it is possible to reach the notion of shape of a basic gesture directly from the Euclidean vector space. As a matter of fact a Euclidean vector space is a particular Quadratic space; that is a vector space with a non-degenerate symmetric bilinear form\footnote{Which, moreover, is positive definite.}. To each Quadratic space it is associated a unique Clifford vector algebra that we denote by the symbol $\mathcal{C}\ell(p,q)$, where $(p,q)\in\mathbb{N}^2$ is the signature of the not-degenerate quadratic form associated to the symmetric bilinear form. Here, we are just interested in the Euclidean plane case: $p=2$ and $q=0$. The positive definite symmetric bilinear form, between vectors $\vec{u}$ and $\vec{v}$ in the two dimensional Euclidean vector space~$\mathbb{E}_2$, is denoted by a dot: $\vec{u}\cdot \vec{v}$. Let $\{\vec{e}_1, \vec{e}_2\}$ be a fixed orthonormal basis for~$\mathbb{E}_2$.  Then, every element $\mathbf{X}\in \mathcal{C}\ell(2,0)$ can be uniquely expressed as
\[
\mathbf{X}=\alpha + x_1 \vec{e}_1 + y \vec{e}_2 + \beta \mathbf{I}\ ,
\]
where $\alpha, x_1,x_2,\beta$ are real numbers, and $\mathbf{I}=\vec{e}_1\vec{e}_2$. With such representation, the Euclidean space $\mathbb{E}_2$ is a vector  subspace of $\mathcal{C}\ell(2,0)$, and the field of scalars $\mathbb{R}$ is a subalgebra of $\mathcal{C}\ell(2,0)$. The associative and distributive (but not necessarily commutative) product in $\mathcal{C}\ell(2,0)$ is uniquely generated by the following simple rule:
\[
\vec{v} \vec{v} = \vec{v}\cdot \vec{v}\ 
\]
for every $\vec{v}\in\mathbb{E}_2\subset \mathcal{C}\ell(2,0)$. \\
The foregoing rule has several important consequences. In particular,

\begin{itemize}
	\item $\displaystyle \frac{1}{2}\big( \vec{u}\vec{v}+\vec{v}\vec{u}\big)=\vec{u}\cdot \vec{v}$, for each $\vec{u}, \vec{v}\in\mathbb{E}_2$;
	\item $\big(\vec{e}_1\big)^2=1=\big(\vec{e}_2\big)^2$, and $\vec{e}_1\vec{e}_2=-\vec{e}_2\vec{e}_1$, $\mathbf{I}^2=-1$;
	\item $\vec{e}_1 \mathbf{I}= \vec{e}_1\vec{e}_1\vec{e}_2= \vec{e}_2\in\mathbb{E}_2$ and, similarly, $\vec{e}_2 \mathbf{I}= -\vec{e}_1\in\mathbb{E}_2$; 
	\item if we define $\displaystyle \vec{u}\wedge \vec{v}=\frac{1}{2}\big( \vec{u}\vec{v}-\vec{v}\vec{u}\big)$, then we have that
		\begin{itemize}
			\item $\displaystyle \vec{u}\wedge \vec{v} 
			= \det\left(
			\begin{array}{cc} u_1 & u_2 \\ v_1 & v_2\end{array}
			\right)\mathbf{I}$, where $\vec{u}= u_1\vec{e}_1+u_2\vec{e}_2$, $\vec{v}= v_1\vec{e}_1+v_2\vec{e}_2$;
			\item $\vec{u}\vec{v}=\vec{u}\cdot \vec{v}+\vec{u}\wedge \vec{v}$.
	  \end{itemize}
\end{itemize}

Then, the associative and distributive product $\mathbf{U}\mathbf{V}$ between  the Clifford numbers $\mathbf{U}=u_0 + \vec{u} + u_{3} \mathbf{I}$ and $\mathbf{V}=v_0 + \vec{v} +v_{3} \mathbf{I}$ in $\mathcal{C}\ell(2,0)$, can be written as follows
\begin{align}
&
u_0 v_0 + \vec{u}\cdot \vec{v}-u_3v_3 +\label{eq:scalar part} \\
+ &
u_0\vec{v} +v_0\vec{u}+v_3 \vec{u}\mathbf{I}+u_3\mathbf{I}\vec{v}+\label{eq:vector part}\\
+ & 
\vec{u}\wedge \vec{v} + u_0v_3\mathbf{I} + u_3 v_0 \mathbf{I}\ ,\label{eq:pseudoscalar part}
\end{align}
where $(\ref{eq:scalar part})\in\mathbb{R}$, $(\ref{eq:vector part})\in\mathbb{E}_2$, and $(\ref{eq:pseudoscalar part})$ is scalar multiple of $\mathbf{I}$. 
Moreover,
\[
\mathbf{U}\cdot \mathbf{V}
=
u_0 v_0 + \vec{u}\cdot \vec{v}+u_3v_3 
\]
defines a positive definite bilinear form on $\mathcal{C}\ell(2,0)$, whose  Euclidean norm\footnote{Note that, if $\mathbf{U}\in\mathbb{R}$, then $\left|\mathbf{U}\right|_{\mathcal{C}\ell(2,0)}=\left|\mathbf{U}\right|_{\mathbb{R}}$, and if $\mathbf{U}\in\mathbb{E}_2$, then $\left|\mathbf{U}\right|_{\mathcal{C}\ell(2,0)}=\left|\mathbf{U}\right|_{\mathbb{E}_2}$.} is
\[
\left|
\mathbf{U}
\right|_{\mathcal{C}\ell(2,0)}
=
\sqrt{(u_0)^2 + \big(|\vec{u}|_{\mathbb{E}_2}\big)^2+ (u_3)^2}\ .
\]
This apparently messy situation hide an algebraic structure that is richer than that of complex numbers, and can encode many different geometric notion of the Euclidean plane within a single algebraic frame. Here, we want to point out just few properties:
\begin{itemize}
	\item every non-zero vector $\vec{v}\in\mathbb{E}_2$ is invertible in $\mathcal{C}\ell(2,0)$ and $\displaystyle (\vec{v})^{-1}=\frac{1}{|\vec{v}|_{\mathbb{E}_2}^2}\vec{v}$
	\item if $\vec{u},\vec{v}\in\mathbb{E}_2$, and $\vec{v}\ne \vec{0}$, then 
	$\displaystyle \vec{u}(\vec{v})^{-1}
	=
	\vec{u}/\vec{v}
	=
	\frac{\vec{u}\cdot\vec{v}}{|\vec{v}|_{\mathbb{E}_2}^2}+\frac{1}{|\vec{v}|_{\mathbb{E}_2}^2}\vec{u}\wedge \vec{u}
	$.
\end{itemize}
Since $\mathbf{i}\in\mathbb{C}$ has the same algebraic properties of $-\mathbf{I}\in \mathcal{C}\ell(2,0)$, we can consider  the shape of a basic gesture $\big(\vec{u},\vec{v}\big)$ ``directly'' as the quotient, in the Clifford algebra $\mathcal{C}\ell(2,0)$, of the two vectors; as a matter of fact 
\[
\mathbb{C}\ni
\frac{\mathbf{u}}{\mathbf{v}}=\frac{rx+sy}{x^2+y^2}-\mathbf{i}\frac{ry-sx}{x^2+y^2}
\longleftrightarrow  \frac{rx+sy}{x^2+y^2}+\mathbf{I}\frac{ry-sx}{x^2+y^2} =\vec{u}/\vec{v}
\in \mathcal{C}\ell(2,0)
\ ,
\]
where $\mathbf{u}=r+\mathbf{i}s$, $\mathbf{v}=x+\mathbf{i}y$, $\vec{u}=r\vec{e}_1+s\vec{e}_2$, and $\vec{v}=x\vec{e}_1+y\vec{e}_2$.
\begin{definition}
The \textbf{shape} of a basic gesture $(\vec{v}_1,\vec{v}_2)$ is the Clifford number 
\[
\vec{v}_1 \big(\vec{v}_2\big)^{-1}
=
\vec{v}_1 \Big/ \vec{v}_2\ \in \mathcal{C}\ell(2,0)\ .
\]
\end{definition}
\begin{definition}
The \textbf{Shape} of a plane gesture $\vec{g}(t)\in\mathbb{E}_2$, is the following multivector-valued function
\[
1
-
\frac{1}{2}
\Big(\vec{g}\ ''(t) \Big/\vec{g}\ '(t)\Big) \in \mathcal{C}\ell(2,0)\ .
\]
\end{definition}
In order to express our previous convergence theorems in terms of Clifford numbers, it suffices to rewrite the Local Shape Distance in terms of scalar products\footnote{As was done in~\cite{2018-!FTL1} for the Javascript implementation of the  algorithm (see appendix~{\bf B.1} to~\cite{2018-!FTL1}).}.

\begin{prop*}
Given two basic gestures $(\vec{u}_1,\vec{u}_2)$ and $(\vec{v}_1,\vec{v}_2)$, then 
\begin{align*}
&
LSD\big((\vec{u}_1,\vec{u}_2),(\vec{v}_1,\vec{v}_2)\big)
= 
\left| \vec{u}_1\Big/ \vec{u}_2\ - \ \vec{v}_1 \Big/ \vec{v}_2 \right|_{\mathcal{C}\ell(2,0)}=\\
=
&
\sqrt{\scriptstyle
\frac{|\vec{u}_1|^2|\vec{v}_2|^2+|\vec{u}_2|^2|\vec{v}_1|^2-2\Big[(\vec{u}_1\cdot\vec{u}_2)(\vec{v}_1\cdot\vec{v}_2)-(\vec{u}_1\cdot\vec{v}_2)(\vec{u}_2\cdot\vec{v}_1)+(\vec{u}_1\cdot\vec{v}_1)(\vec{u}_2\cdot\vec{v}_2)\big]}
{|\vec{u}_2|^2|\vec{v}_2|^2}
}\ ,
\end{align*}
where $|\cdot|=|\cdot |_{\mathbb{E}_2}$.
\end{prop*}


\end{document}